\documentclass[reqno]{amsart}
\usepackage{amsmath}
\usepackage{amsthm}
\usepackage{amssymb}
\usepackage{epsfig}
\usepackage{psfrag}
\usepackage{graphicx}
\usepackage{caption2}
\usepackage{textcomp}
\usepackage{pgf}

\newcommand{\Rz}{\mathbb{R}}
\newcommand{\Nz}{\mathbb{N}}

\newcommand{\disp}{\displaystyle}

\newcommand{\ove}{\overline}
\newcommand{\epsi}{\varepsilon}

\newcommand{\e}{\text{\rm e}}
\newcommand{\define}{:=}

\newcommand{\UU}{{\mathcal U}}
\newcommand{\VV}{{\mathcal V}}
\newcommand{\KK}{{\mathcal K}}
\newcommand{\media}{\disp -\!\!\!\!\!\!\int}

  \newtheorem{theorem}{Theorem}[section]

 \newtheorem{lemma}[theorem]{Lemma}
  \newtheorem{conjecture}[theorem]{Conjecture}
\begin{document}

\title[De Giorgi conjecture]{The De Giorgi conjecture on \\ elliptic regularization} 

\author{Ulisse Stefanelli}
\address{IMATI - CNR, 
v. Ferrata 1, I-27100 Pavia, Italy.}
\email{ulisse.stefanelli\,@\,imati.cnr.it} 
\urladdr{http://www.imati.cnr.it/ulisse/}
\thanks{U.S. is partially supported by FP7-IDEAS-ERC-StG Grant \#200497 (BioSMA)}
\keywords{Semilinear wave equation, Elliptic regularization}

\begin{abstract}
We prove a conjecture by De Giorgi on the elliptic regularization of semilinear wave equations in the finite-time case.
\end{abstract}

\subjclass{35L71} 
\maketitle

\pagestyle{myheadings}

                                %
                                %
                                %

\section{Introduction}
\setcounter{equation}{0}

In \cite{Degi95CCEP} {\sc De Giorgi} proposed  the following conjecture. 

\begin{conjecture}\label{1.1}
Let be given $u^0, \ u^1 \in C^\infty_0(\Rz^d)$, $p=2k$, $k>1$ integer. For every $\epsi >0$ let $u^\epsi$ be a minimizer of the functional for $T=\infty$ 
\begin{equation*}
  \label{I}
  I_\epsi(u)\define \int_{\Rz^d \times (0,T)} \e^{-t/\epsi}\left(  |u_{tt}|^2 +  {\epsi^{-2}}|\nabla u|^2 +  {\epsi^{-2}} |u|^p  \right) dx\, dt 
\end{equation*}
in the class of all functions $u$ satisfying the initial conditions $u^\epsi(x,0)=u^0(x)$, $u^\epsi_t(x,0)=u^1(x)$. Then, there exists $\lim_{\epsi \to 0} u^\epsi(x,t) =  u(x,t)$, and satisfies the equation
\begin{equation}
  \label{PDE}
  u_{tt} - \Delta u + \frac{p}{2}|u|^{p-2}u = 0.
\end{equation}
\end{conjecture}

The interest in this conjecture resides in the possibility of connecting the {\it difficult} semilinear wave equation \eqref{PDE} with a comparably {\it easier} problem: the constrained minimization of the uniformly convex functional $I_\epsi$. Note for instance that no uniqueness for the semilinear wave equation \eqref{PDE} is available for large $p$ whereas the functional $I_\epsi$ always admits a unique minimizer.

The reference to {\it elliptic regularization} in this context is related to the fact that the Euler-Lagrange equation for the functional  $I_\epsi$ can be formally computed as 
\begin{equation}
  \label{EL}
  \epsi^2 u_{tttt}^\epsi - 2 \epsi u_{ttt}^\epsi + u^\epsi_{tt} - \Delta u^\epsi + \frac{p}{2}  |u^\epsi|^{p-2}u^\epsi = 0.
\end{equation}
Namely, minimizing $I_\epsi$ consists in addressing a fourth-order elliptic regularization in time of the semilinear wave equation \eqref{PDE}. 

 The statement of Conjecture \ref{1.1} does not specify the convergence notion for $u^\epsi \to u$ nor the solution notion for equation \eqref{PDE}. We shall tacitly assume in the following that the convergence is (at least) almost everywhere in space and time,  that convergence holds (at least) for subsequences, and that equation \eqref{PDE} will be solved (at least) in the distributional sense. Along with these provisions, the main result of the paper is the following.

\begin{theorem}\label{main}
 Let $T < \infty$. Then,  Conjecture \emph{\ref{1.1}} holds true.
\end{theorem}

\begin{figure}[h]
 \centering
\includegraphics[width=0.8\linewidth]{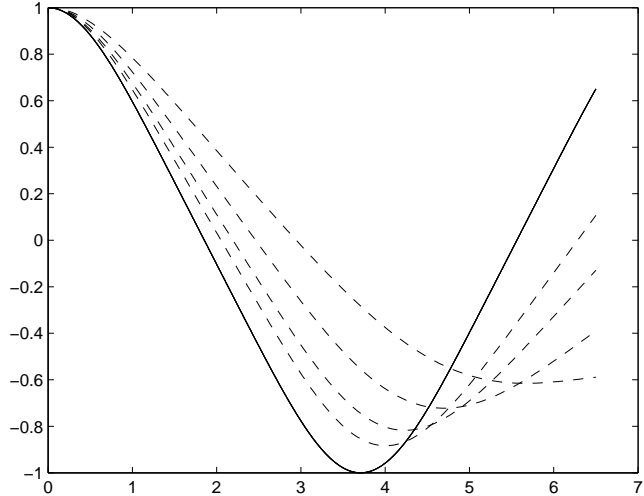}
 \caption{Conjecture \ref{1.1} in the scalar case $d=0$, for $p=4$, $u^0=1$, and $u^1=0$. The dashed lines are the minimizers of $I_\epsi$ for $\epsi=0.4$ (top at $t=0+$), $0.2$, $0.1$, $0.05$ (bottom), whereas the solid line is the unique solution of $u''+u^3=0$, along with the given boundary conditions.}
\label{fig}
\end{figure}

In the scalar case $d=0$, Conjecture \ref{1.1} is illustrated in Figure \ref{fig}.

Note that the possibility of considering $T =\infty$ seems to be nontrivial and it is currently out of reach for the present analysis.  Moreover, although not explicitly mentioned in the original text of Conjecture \ref{1.1}, one may quite naturally wonder if the limit trajectory $u$ attains the initial conditions $u(0)=u^0$ and $u_t(0)=u^1$ as well. In this regard, we shall remark that our analysis is still partial as the only condition $u(0)=u^0$ is ascertained in the following.

On the other hand, we are in the position of proving some generalization of Conjecture \ref{1.1} as well. Indeed, assumptions on the initial data $u^0$ and $u^1$ can be weakened and convergence is proved to take place in a suitable topology (and, correspondingly, to a stronger solution notion).  More significantly, we can allow for more general nonlinearities, in the spirit of  \cite{Degi95CCEP} . In particular, we prove Conjecture \ref{1.1} for all $p >2$. 

The natural energy estimate for semilinear wave equation \eqref{PDE} entails a pointwise-in-time bound on $u$ in terms of data.
The core of the proof of Conjecture \ref{1.1} relies in checking that, even in the elliptic-regularized situation of the Euler-Lagrange equation \eqref{EL}, we are in the position of establishing a corresponding {\it integral} energy estimate for $u^\epsi$ independently of $\epsi$. Namely, we have the following.

\begin{lemma}[Energy estimate]\label{ee} Let $u^\epsi$ minimize $I_\epsi$ for $\epsi<1/2$ with  $u^\epsi(0)=u^0$ and $u^\epsi_t(0)=u^1$. Then,
  \begin{equation}
    \label{ap}
    \int_{\Rz^d \times (0,T)} \left( |u_{t}^\epsi|^2 + |\nabla u^\epsi|^2 + |u^\epsi|^p\right)  dx \, dt  \leq C
  \end{equation}
where $C$ depends on $u^0$, $u^1$, and $T$ but not on $\epsi$.
\end{lemma}

The proof of this Lemma will be detailed in Section \ref{apri}. Once the energy estimate \eqref{ap} is established, Conjecture \ref{1.1} follows by passing to the limit in (a suitable variational version of) the Euler-Lagrange equation \eqref{EL}. This limiting  procedure is discussed in Section \ref{limit}. Eventually, extensions and comments are collected in Subsection \ref{comments}.

This is, to our knowledge, the first result in the direction of elliptic regularization of nonlinear hyperbolic equations. Elliptic regularization has been already considered in the frame of parabolic problems and results in the linear case can be found in {\sc Lions \& Magenes} \cite{Lions-Magenes1F}. The first application of this variational perspective to a nonlinear dissipative problem is due to {\sc Ilmanen} \cite{Ilmanen94} in the context of mean curvature flow. Results and applications to rate-independent dissipative systems have been presented by {\sc Mielke \& Ortiz} \cite{Mielke-Ortiz06} and then extended and coupled with time-discretization in \cite{ms2}. In the case of gradient flows, two relaxation and scaling examples are provided by {\sc Conti \& Ortiz} \cite{Conti-Ortiz08} whereas the general abstract theory is addressed in \cite{ms3}. Some application to mean curvature flow of cartesian surfaces is in \cite{spst} and the elliptic regularization of doubly nonlinear parabolic equations is discussed in \cite{akst,akst2}. Finally, a similar functional approach (with $\epsi$ fixed though) has been considered by {\sc Lucia, Muratov, \& Novaga} \cite{Lucia-et-al08,Muratov-Novaga08,Muratov-Novaga08b} in connection with travelling waves in reaction-diffusion-advection problems.

On the other hand, the literature on the semilinear wave equation \eqref{PDE} is vast and it is clearly beyond the purposes of this note to provide a comprehensive review. The reader can consider the monographs by {\sc Lions} \cite[Chap. 1]{LionsNL} and by {\sc Shatah \& Struwe} \cite[Subsect. 8.3.1]{Shatah-Struwe98} for a collection of results, references, and historical remarks.

\section{Energy estimate}\label{apri}
\setcounter{equation}{0}

This section brings to the proof of Lemma \ref{ee}. We start by presenting a formal argument in Subsection \ref{formal}. Then, the rigorous proof is achieved by means of a time-discretization procedure and is detailed in the Subsections 2.2-2.5.

In the following $(\cdot, \cdot)$ stands for the usual scalar product in $L^2(\Rz^d)$ and $|\cdot|$ is used both for the corresponding norm and the modulus in $\Rz$ and $\Rz^d$. Moreover, given
$$U:=\{u \in L^p(\Rz^d) \ : \ \nabla u \in L^2(\Rz^d)\}$$
(note that $U \not \subset L^2(\Rz^d)$, take $x \mapsto  (1+ |x|)^{-d/2}$),
for all $u, \, v \in U$ we will use the notation
\begin{align}
\phi(u) &\define  \frac12 \int_{\Rz^d}\left( |\nabla u|^2 + |u|^p\right),  \nonumber\\
(Au,v)&\define\int_{\Rz^d} \left(\nabla u \cdot \nabla v + \frac{p}{2}|u|^{p-2}u\,v \right).  \nonumber
\end{align}
In particular, let us observe that $(p'\define p/(p-1)<2<p)$
\begin{align}
 (Au,v)&\leq \int_{\Rz^d}\left(\frac{\sigma^2}{2} |\nabla u|^2 + \frac{1}{2\sigma^2}|\nabla v|^2 + \frac{\sigma^{p'}}{p'}  \frac{p}{2}  \Big||u|^{p-2}u\Big|^{p'} + \frac{1}{p\sigma^p}  \frac{p}{2}  |v|^p\right)\nonumber\\
&\leq \sigma^{p'} (p-1)\phi(u) + \frac{1}{\sigma^p} \phi(v) \quad  \text{for all} \ \  \sigma \in (0,1).\label{somma_phi}
\end{align}

Define the space
\begin{align}
 \UU\define H^2(0,T;L^2(\Rz^d))\cap L^2(0,T;H^1(\Rz^d))\cap L^p(\Rz^d\times (0,T))\nonumber
\end{align}
and let $u^\epsi$ minimize $I_\epsi$ over the non-empty, convex, and closed set
$$\KK(u^0,u^1)\define \{u \in \UU \ : \ u(0)=u^0, \ u_t(0)=u^1\}.$$
By considering $r \mapsto I_\epsi(u^\epsi + r v)$ for $v \in K(0,0)$ we obtain that
\begin{equation}
  \label{var}
  0= \int_0^T \e^{-t/\epsi}\left((u^\epsi_{tt},v_{tt}) + \frac{1}{\epsi^2}(Au^\epsi,v) \right) \quad \forall v \in \KK(0,0).
\end{equation}
For the sake of notational simplicity, we shall drop the superscript $\epsi$ from $u^\epsi$ in the remainder of this section.

\subsection{A formal estimate}\label{formal}
Let us present here some heuristics for Lemma \ref{ee}. By assuming smoothness, we get from \eqref{var} that $u$ solves 
\begin{align}
 0&= \int_0^T \left( \partial_{tt}(\e^{-t/\epsi} u_{tt}) + \frac{1}{\epsi^2}\e^{-t/\epsi} \left(-\Delta u +  \frac{p}{2}  |u|^{p-2}u\right) ,v\right) \nonumber\\
&+ \big(\e^{-T/\epsi}u_{tt}(T),v_t(T)\big) - \big(\partial_t(\e^{-t/\epsi}u_{tt})(T), v(T) \big)
\quad \forall v \in \KK(0,0).\label{formal0}
\end{align}
Namely, the Euler-Lagrange equation \eqref{EL} holds along with the {\it final} homogeneous Neumann boundary conditions
\begin{equation}
  \label{fbc}
  u_{tt}(T)= u_{ttt}(T) =0.
\end{equation}

Test now \eqref{EL} by the function $t \mapsto (1 +T - t) (u_t(t) - u^1) $ and take the integral on $(0,T)$. By recalling that 
\begin{equation}\label{gg}\forall g \in L^1(0,T), \quad \int_0^T (1 + T - t)g(t) dt = \int_0^T g(t)dt + \int_0^T \left(\int_0^t g(s)ds\right) dt
\end{equation}
we obtain the equality 
\begin{align}
 &\left(\frac12 - \epsi\right)|u_t(T) - u^1|^2 + \frac12 \int_0^T |u_t - u^1|^2 + \phi(u(T)) + \int_0^T \phi(u) \nonumber\\
&+{} \frac{(1+T)\epsi^2}{2}|u_{tt}(0)|^2 + \left( 2\epsi - \frac32 \epsi^2 \right)\int_0^T |u_{tt}|^2 + 2\epsi \int_0^T \int_0^t|u_{tt}|^2 \nonumber\\
&= -{}\epsi^2 (u_{ttt}(T),u_t(T)-u^1) + \frac{\epsi^2}{2}|u_{tt}(T)|^2 + (2\epsi - \epsi^2)(u_{tt}(T), u_t(T)-u^1) \nonumber\\
&+ (1+T)\phi(u^0)+ \int_0^T (Au,u^1)+ \int_0^T\int_0^t (Au,u^1).\label{for}
\end{align}
Note that this computation amounts in testing \eqref{EL} by $t \mapsto u_t(t) -u^1$, integrating in time, and summing the result with the same relation integrated in time once again (see \eqref{gg}). Namely, we summarize this procedure as follows
\begin{equation}\label{spiega}
 \int_0^T \eqref{EL}\big|_{v=u_t - u^1}   \ + \ \int_0^T \left(\int_0^t \eqref{EL}\big|_{v=u_t - u^1}\right).
\end{equation}

Owing to the final boundary conditions \eqref{fbc}, the right-hand side of \eqref{for} can be controlled as follows
\begin{gather}
(1+T)\phi(u^0) + \int_0^T (Au,u^1)+ \int_0^T\int_0^t (Au,u^1)\nonumber\\
\stackrel{\eqref{somma_phi}}{\leq} \frac12 \int_0^T \phi(u) + C\Big( \phi(u^0) + \phi(u^1)\Big)\label{rhs}
\end{gather}
where $C$ depends on $T$. Hence, for $\epsi<1/2$, one has that 
$$ \int_0^T \Big( |u_t|^2 +  \phi(u) \Big) \leq C.$$
Namely, the energy estimate \eqref{ap} holds.

\subsection{Time-discretization}

The argument of Subsection \ref{formal} is formal. In particular, the final boundary conditions \eqref{fbc} make no sense as $u\in \UU$ only. 

In order to provide a rigorous proof of Lemma \ref{ee} we will proceed by time-discretization. To this aim let the time step $\tau\define T/n \, (n \in \Nz)$ be fixed, consider the space
$$\UU_\tau\define \big\{(u_0,\dots,u_n) \in (L^2(\Rz^d))^{n+1}\ : \ (u_2,\dots,u_{n-2}) \in U^{n-3}\big\},$$
and define the functional $I_{\epsi\tau}: \UU_\tau \to [0,\infty)$  by
\begin{equation*}
  \label{Itau}
  I_{\epsi\tau}(u_0,\dots,u_n) \define \sum_{i=2}^n \tau \rho_{\tau,i} \frac12 |\delta^2u_i|^2 +  \sum_{i=2}^{n-2} \frac{\tau}{\epsi^2} \rho_{\tau,i+2} \phi(u_i).
\end{equation*}
Given the vector $(w_m,\dots,w_n)$, in the latter we have used the notation $\delta w$ for its {\it discrete derivative} $\delta w_i\define (w_i - w_{i-1})/\tau$ for $i\geq m+1$ and $\delta^2w\define \delta(\delta w)$, $\delta^3w \define \delta(\delta^2w)$ and so on. Moreover, we have used the {\it weights}
\begin{equation}
  \label{rho}
   \rho_{\tau,i} \define \left( \disp \frac{\epsi}{\epsi + \tau}\right)^i\quad i=0,1,\dots,n.
\end{equation}
These weights are nothing but the discrete version of the exponentially decaying weight $t \mapsto \exp(-t/\epsi)$ for we have that $\delta\rho_{\tau,i} + \rho_{\tau,i}/\epsi =0$. Namely, $\rho_{\tau,i}$ is the solution of the backward Euler scheme applied to $\rho'+\rho/\epsi=0$. 

The functional $I_{\epsi\tau}$ represents a discrete version of the original time-continuous functional $I_\epsi$. We shall drop the subscript $\tau$ from $\rho_{\tau,i}$ in the remainder of this section for notational simplicity.

For all vectors  $(w_0,\dots,w_n)$, we indicate by $\ove w_\tau$ and $w_\tau$ its backward constant and piecewise affine interpolants on the partition, respectively. Namely,
\begin{gather}
  \ove w_\tau(0)=w_\tau(0)\define w_0, \quad \ove w_\tau(t)\define w_i,\quad w_\tau(t)\define \alpha_i(t)w_i + (1-\alpha_i(t))w_{i-1}\nonumber  \\
 \quad \text{for} \ \ t \in ((i-1)\tau,i\tau], \ i=1,\dots,n\nonumber
\end{gather}
where we have used the auxiliary functions 
$$  \alpha_i(t)\define (t- (i-1)\tau)/\tau \quad \text{for} \ \ t \in ((i-1)\tau,i\tau], \ i=1,\dots,n.$$
In particular, $\partial_t w_\tau = \ove{\delta w}_\tau$ almost everywhere. 

The strategy of the proof of Lemma \ref{ee} consists in establishing the energy estimate \eqref{ap} at the time-discrete level (Subsection \ref{stab}) and then pass to the limit as $\tau \to 0$ (Subsection \ref{passage}). As a preparatory step, we present the Euler-Lagrange equation for the time-discrete functional $I_{\epsi\tau}$ in Subsection \ref{eesec}.

\subsection{Time-discrete Euler-Lagrange equation}\label{eesec}
Let us define the convex set
$$\KK_\tau(u^0,u^1)\define \big\{ (u_0,\dots,u_n) \in \UU_\tau  \ : \  u_0=u^0, \ \delta u_1=u^1\big\}.$$
The time-discrete functional $I_{\epsi\tau}$ is clearly convex. Moreover, $I_{\epsi\tau}$ turns out to be coercive on $\KK_\tau(u^0,u^1)$. Indeed, the coercivity of $I_{\epsi\tau}$ in $U^{n-3}$ with respect to the components $(u_2,\dots,u_{n-2})$ is immediate. As for the coercivity in $(L^2(\Rz^d))^{n+1}$ we preliminarily observe that, for all $(u_0,\dots,u_n) \in \KK_\tau$,
\begin{equation}
\sum_{k=2}^n\tau |u_i|^2 \leq C\left( |u^0|^2 + |u^1|^2 + n \sum_{k=2}^n\tau |\delta^2 u_k|^2\right)\quad  i=2, \dots,n \label{poin}
\end{equation}
where $C$ depends on $T$. Hence, we have that
\begin{align}
I_{\epsi,\tau}(u_0,\dots,u_n) \geq \sum_{k=2}^n \tau \rho_k \frac12 |\delta^2 u_k|^2 
\stackrel{\eqref{poin}}{\geq} \frac{\rho_n}{2n}\left( \frac{1}{C} \sum_{k=2}^n \tau |u_k|^2 - |u^0|^2 - |u^1|^2\right)\nonumber
\end{align}
where $C$ is the constant in \eqref{poin}. In particular, the functional $I_{\epsi\tau}$ controls the norm in $(L^2(\Rz^d))^{n+1}$
 (up to constants depending on $T$, $\epsi$, and $\tau$). Moreover, the quadratic character of the first term in $I_{\epsi\tau}$ and the latter computation ensure that the functional $I_{\epsi\tau}$ is uniformly convex in $(L^2(\Rz^d))^{n+1}$ along $\KK_\tau(u^0,u^1)$. Hence, $I_{\epsi\tau}$ admits a unique minimizer $(u_0,\dots,u_n)$ in $ \KK_\tau(u^0,u^1) $ and we directly compute (see equation \eqref{var}) that 
\begin{gather}
  0 = \sum_{i=2}^n \tau\rho_{i} (\delta^2u_i, \delta^2v_i) + \sum_{i=2}^{n-2}\frac{\tau}{\epsi^2} \rho_{i+2} (Au_i,v_i) \quad \forall (v_0,\dots,v_n)\in \KK_\tau(0,0). \label{vartau}
\end{gather}

Let us now proceed along the same lines of \eqref{formal0}. First of all, we sum-by-parts and obtain that
\begin{align}
 &\sum_{i=2}^n \tau\rho_{i} (\delta^2u_i, \delta^2v_i) = \sum_{i=2}^n \rho_{i} (\delta^2u_i, \delta v_i) - \sum_{i=2}^n \rho_{i} (\delta^2u_i, \delta v_{i-1})\nonumber\\
&= \rho_{n}(\delta^2 u_n,\delta v_n) - \rho_{2}(\delta^2 u_2, \delta v_1) - \sum_{i=2}^{n-1} \tau (\delta(\rho\delta^2u)_{i+1}, \delta v_i)\nonumber\\
&\stackrel{\delta v_1=0}{=} \rho_{n}(\delta^2 u_n,\delta v_n) - \sum_{i=2}^{n-1} (\delta(\rho\delta^2u)_{i+1},  v_i) +  \sum_{i=2}^{n-1} (\delta(\rho\delta^2u)_{i+1},  v_{i-1})\nonumber\\
&= \rho_{n}(\delta^2 u_n,\delta v_n)  - (\delta (\rho \delta^2 u)_n, v_{n-1}) + (\delta (\rho \delta^2 u)_3,v_1) +  \sum_{i=2}^{n-2} \tau (\delta^2(\rho\delta^2u)_{i+2},  v_i)\nonumber\\
&\stackrel{v_1=0}{=}\rho_{n}(\delta^2 u_n,\delta v_n)  - (\delta (\rho \delta^2 u)_n, v_{n-1}) +  \sum_{i=2}^{n-2} \tau (\delta^2(\rho\delta^2u)_{i+2},  v_i).\nonumber
\end{align}
Next, by means of definition \eqref{rho} (and some tedious computation) we check that
\begin{align}
  \delta^2(\rho\delta^2u)_{i+2} = \rho_{i+2}\left( \delta^4u_{i+2} - \frac{2}{\epsi}\delta^3u_{i+1} + \frac{1}{\epsi^2}\delta^2 u_i\right),\nonumber
\end{align}
and rewrite relation \eqref{vartau} in the equivalent form
\begin{align}
  0 &= \rho_{n}(\delta^2 u_n,\delta v_n)  - (\delta (\rho \delta^2 u)_n, v_{n-1}) \nonumber\\
&+  \sum_{i=2}^{n-2}\frac{\tau}{\epsi^2} \rho_{i+2}\left( (\epsi^2\delta^4u_{i+2} -  2\epsi\delta^3u_{i+1} + \delta^2 u_i,v_i) + (Au_i ,v_i)\right)  \nonumber\\
 &\text{for all} \ \ (v_0,\dots,v_n)\in \KK_\tau(0,0). \label{vartau2}
\end{align}

We have hence proved that the minimizer $(u_0, \dots, u_n)$ of $I_{\epsi\tau}$ solves 
\begin{gather}
(\epsi^2\delta^4u_{i+2} -  2\epsi\delta^3u_{i+1} + \delta^2 u_i,v) + ( A u_i, v)=0 \quad
 \forall v \in U, \ i= 2,\dots,n-2,\label{ELtau}\\
u_0=u^0, \quad \delta u_1=u^1, \quad \delta^2 u_n = \delta^3 u_n =0,\label{fbctau}  
\end{gather}
which is nothing but the time-discrete analogue of the Euler-Lagrange equation \eqref{EL} and of the final conditions \eqref{fbc}. Note incidentally that the final boundary conditions in \eqref{fbctau} imply, in particular, that $\delta^2 u_{n-1}=0$.

\subsection{Energy estimate at the time-discrete level}\label{stab}
The strategy of the proof of the energy estimate at the time-discrete level is exactly the same as that of Section \ref{formal}. In particular, we present a time-discrete version of estimate \eqref{spiega} by using the time-discrete Euler-Lagrange equation \eqref{ELtau} instead of \eqref{EL}. Namely, we perform the following
\begin{equation}\label{spiegatau}
 \sum_{k=2}^{n-2} \tau\eqref{ELtau}\big|_{v=\delta u_k - u^1}  \ + \  \sum_{k=2}^{n-2}\tau \left( \sum_{i=2}^{k}\tau \eqref{ELtau}\big|_{v=\delta u_i - u^1}\right),
\end{equation}
(note that $\delta u_k - u^1 \in U$ for $k=2, \dots, n-2$).

We shall start by the first term in \eqref{spiegatau}. Hence, we test relation \eqref{ELtau} on $v =\tau( \delta u_i - u^1)$ and sum for $i=2,\dots,k\leq n-2$ in order to get that 
\begin{align}
&\epsi^2 \sum_{i=2}^k \tau (\delta^4 u_{i+2},  \delta u_i - u^1) -2\epsi \sum_{i=2}^k\tau (\delta^3 u_{i+1},  \delta u_i - u^1) \nonumber\\
&{}+ \sum_{i=2}^k \tau (\delta^2 u_{i},  \delta u_i - u^1) +\sum_{i=2}^k \tau (A u_i, \delta u_i - u^1)=0. \label{complete}  
\end{align}

Let us now treat separately the terms in the above left-hand side. The fourth-order-in-time term can be handled as follows.
\begin{align}
 &\epsi^2 \sum_{i=2}^k \tau (\delta^4 u_{i+2},  \delta u_i - u^1) = \epsi^2 \sum_{i=2}^k (\delta^3 u_{i+2} - \delta^3 u_{i+1}, \delta u_i - u^1)\nonumber\\
&=\epsi^2 (\delta^3 u_{k+2}, \delta u_k - u^1) - \epsi^2 (\delta^3 u_3, \delta u_2 - u^1) - \epsi^2 \sum_{i=3}^k \tau (\delta^3 u_{i+1}, \delta^2 u_i) \nonumber\\
&\stackrel{\delta u_1 =u^1}{=} \epsi^2 (\delta^3 u_{k+2}, \delta u_k - u^1) - \epsi^2 \sum_{i=2}^k \tau (\delta^3 u_{i+1}, \delta^2 u_i) \nonumber\\
&= \epsi^2 (\delta^3 u_{k+2}, \delta u_k - u^1)  -\frac{\epsi^2}{2}|\delta^2 u_{k+1}|^2 \nonumber\\
&+ \frac{\epsi^2}{2}|\delta^2 u_2|^2 + \epsi^2 \sum_{i=2}^k |\delta^2 u_{i+1} - \delta^2 u_i|^2.\label{o4}
\end{align}

Next, we treat the third-order-in-time term of \eqref{complete} as
\begin{align}
 &-2\epsi \sum_{i=1}^k\tau (\delta^3 u_{i+1},  \delta u_i - u^1)   =  -2 \epsi \sum_{i=2}^k (\delta^2 u_{i+1} - \delta^2 u_i, \delta u_i - u^1) \nonumber\\
&= -2\epsi (\delta^2 u_{k+1}, \delta u_k - u^1) + 2\epsi (\delta^2 u_2, \delta u_2 - u^1) + 2 \epsi \sum_{i=3}^k \tau |\delta^2 u_i|^2\nonumber\\
&\stackrel{\delta u_1 =u^1}{=} -2\epsi (\delta^2 u_{k+1}, \delta u_k - u^1) + 2\epsi \sum_{i=2}^k \tau |\delta^2 u_i|^2.\label{o3}
\end{align}

As for the remaining terms in \eqref{complete} we compute 
\begin{align}
& \sum_{i=2}^k\tau (\delta^2 u_i, \delta u_i - u^1) = \frac12 |\delta u_k - u^1|^2 +\frac12 \sum_{i=2}^k |\delta u_i - \delta u_{i-1}|^2,\label{o2}\\
& \sum_{i=2}^k\tau (A u_i, \delta u_i - u^1) \stackrel{\delta u_1=u^1}{=}\sum_{i=1}^k\tau (A u_i, \delta u_i - u^1) \nonumber\\
& \geq \phi(u_k) - \phi(u^0) - \sum_{i=1}^k\tau (A u_i, u^1). \label{o0}  
\end{align}

By recollecting the computations \eqref{o4}-\eqref{o0} into equation \eqref{complete} we deduce that
\begin{align}
&\epsi^2 (\delta^3 u_{k+2}, \delta u_k - u^1)  -\frac{\epsi^2}{2}|\delta^2 u_{k+1}|^2 + \frac{\epsi^2}{2}|\delta^2 u_2|^2  -2\epsi (\delta^2 u_{k+1}, \delta u_k - u^1)  \nonumber\\
& + 2\epsi \sum_{i=2}^k \tau |\delta^2 u_i|^2+\frac12 |\delta u_k - u^1|^2 +  \phi(u_k) 
\leq  \phi(u^0) + \sum_{i=1}^k\tau (A u_i, u^1).\label{complete2}
\end{align}

Let us now move to the consideration of the second term in \eqref{spiegatau}. We multiply \eqref{complete2} by $\tau$ and take the sum for $k=2,\dots,n-2$ getting
\begin{align}
&\epsi^2 \sum_{k=2}^{n-2}\tau(\delta^3 u_{k+2}, \delta u_k - u^1)  -\frac{\epsi^2}{2}\sum_{k=2}^{n-2}\tau |\delta^2 u_{k+1}|^2 + \frac{\epsi^2}{2}(n-3)\tau|\delta^2 u_2|^2  \nonumber\\
&-2\epsi\sum_{k=2}^{n-2}\tau (\delta^2 u_{k+1}, \delta u_k - u^1)  + 2\epsi \sum_{k=2}^{n-2}\tau\sum_{i=2}^k \tau |\delta^2 u_i|^2\nonumber\\
&
+\frac12\sum_{k=2}^{n-2}\tau |\delta u_k - u^1|^2 +  \sum_{k=2}^{n-2}\tau \phi(u_k)  \nonumber\\
&
\leq (n-3)\tau \phi(u^0) + \sum_{k=2}^{n-2}\tau\sum_{i=1}^k\tau (A u_i, u^1).\label{complete3}
\end{align}

Before going on, we prepare some computations in order to handle some terms in the above left-hand side.
By summing by parts, we have that 
\begin{align}
 &\epsi^2 \sum_{k=2}^{n-2}\tau(\delta^3 u_{k+2}, \delta u_k - u^1)  = \epsi^2 \sum_{k=2}^{n-2} (\delta^2 u_{k+2} - \delta^2 u_{k+1}, \delta u_k - u^1)\nonumber\\
& = \epsi^2 (\delta^2 u_n, \delta u_{n-2} - u^1) - \epsi^2 (\delta^2 u_3, \delta u_2 - u^1) - \epsi^2 \sum_{k=3}^{n-2} \tau (\delta^2 u_{k+1}, \delta^2 u_k)\nonumber\\
&\stackrel{\eqref{fbctau}}{=} -  \epsi^2 \sum_{k=2}^{n-2} \tau (\delta^2 u_{k+1}, \delta^2 u_k).\label{poi}
\end{align}

Moreover, we can also compute that
\begin{align}
& - 2\epsi \sum_{k=2}^{n-2}\tau (\delta^2 u_{k+1}, \delta u_k - u^1) =  2\epsi \sum_{k=2}^{n-2} (\delta u_{k} - \delta u_{k+1}, \delta u_k - u^1)\nonumber\\
&= \epsi |\delta u_2 - u^1|^2 + \epsi \sum_{k=2}^{n-3}|\delta u_{k+1} - \delta u_k|^2 - \epsi |\delta u_{n-2} - u^1|^2 \nonumber\\
&+ 2\epsi (\delta u_{n-2} - \delta u_{n-1}, \delta u_{n-2} - u^1)\nonumber\\
&= \epsi |\delta u_2 - u^1|^2 + \epsi \sum_{k=2}^{n-3}|\delta u_{k+1} - \delta u_k|^2 - \epsi|\delta u_{n-2} - u^1|^2 \nonumber\\
&- 2\epsi\tau (\delta^2 u_{n-1},\delta u_{n-2} - u^1)\nonumber\\
&\stackrel{\eqref{fbctau}}{=}\epsi |\delta u_2 - u^1|^2 + \epsi \sum_{k=2}^{n-3}|\delta u_{k+1} - \delta u_k|^2 - \epsi|\delta u_{n-2} - u^1|^2 .\label{poi2}
\end{align}

Finally, let us observe that 
\begin{align}
& 2\epsi \sum_{k=2}^{n-2}\tau |\delta^2 u_k|^2  - \epsi^2 \sum_{k=2}^{n-2}\tau (\delta^2 u_{k+1},\delta^2 u_k) - \frac{\epsi^2}{2} \sum_{k=2}^{n-2}\tau |\delta^2 u_{k+1}|^2 \nonumber\\
&\geq 2\epsi \sum_{k=2}^{n-2}\tau |\delta^2 u_k|^2  - \frac{\epsi^2}{2} \sum_{k=2}^{n-2}\tau |\delta^2 u_{k+1}|^2 - \frac{\epsi^2}{2} \sum_{k=2}^{n-2}\tau |\delta^2 u_{k}|^2   - \frac{\epsi^2}{2} \sum_{k=2}^{n-2}\tau |\delta^2 u_{k+1}|^2 \nonumber\\
& \geq \left(2\epsi - \frac{3}{2}\epsi^2 \right)\sum_{k=2}^{n-2}\tau |\delta^2 u_k|^2 - \epsi^2 \tau |\delta^2 u_{n-1}|^2\nonumber\\
&\stackrel{\eqref{fbctau}}{=}\left(2\epsi - \frac{3}{2}\epsi^2 \right)\sum_{k=2}^{n-2}\tau |\delta^2 u_k|^2.
  \label{poi3}
\end{align}

Let us now write estimate \eqref{complete2} by choosing $k=n-2$ and taking advantage of the final boundary conditions \eqref{fbctau} as
\begin{align}
& \frac{\epsi^2}{2}|\delta^2 u_2|^2 + 2\epsi \sum_{k=2}^{n-2} \tau |\delta^2 u_k|^2+\frac12 |\delta u_{n-2} - u^1|^2 +  \phi(u_{n-2}) 
\leq  \phi(u^0) + \sum_{k=1}^{n-2}\tau (A u_k, u^1).\nonumber
\end{align}
By taking the sum of the latter and \eqref{complete3} and using equality \eqref{poi} and estimates \eqref{poi2}-\eqref{poi3}, we obtain that
\begin{align}
&\left(\frac12 - \epsi \right) |\delta u_{n-2} - u^1|^2 + \frac12 \sum_{k=2}^{n-2}\tau |\delta u_k - u^1|^2 + \phi(u_{n-2}) + \sum_{k=2}^{n-2}\tau \phi(u_k)\nonumber\\
& \frac{\epsi^2}{2}\big( 1 + (n-3)\tau\big) |\delta^2 u_2|^2 + \left(2\epsi - \frac{3}{2}\epsi^2 \right)\sum_{k=2}^{n-2}\tau |\delta^2 u_k|^2  + 2\epsi \sum_{k=2}^{n-2}\tau \sum_{i=2}^k \tau |\delta^2 u_i|^2\nonumber\\
&\leq  \big( 1 + (n-3)\tau\big) \phi(u^0)  + \sum_{k=1}^{n-2}\tau (A u_k, u^1)+ \sum_{k=2}^{n-2}\tau\sum_{i=1}^k\tau (A u_i, u^1).\nonumber
\end{align}
The latter is nothing but the discrete analogue of the former (and formal) estimate \eqref{for}. Similarly to \eqref{rhs}, the above right-hand side can be bounded as follows
\begin{align}
&\big( 1 + (n-3)\tau\big) \phi(u^0)  + \sum_{i=1}^k\tau (A u_i, u^1)+ \sum_{k=2}^{n-2}\tau\sum_{i=1}^k\tau (A u_i, u^1)\nonumber\\
  &\stackrel{\eqref{somma_phi}}{\leq} \frac{1}{2}  \sum_{i=2}^k \tau \phi(u_k) + C\Big(\phi(u^0) + \phi(u^1)\Big).\nonumber
\end{align}
Eventually, we have proved the following.

\begin{lemma}[Discrete energy estimate] Let $(u_0,\dots,u_n)$ minimize $I_{\epsi\tau}$ for $\epsi <1/2$ with  $u_0=u^0$ and $\delta u_1=u^1$. Then,
  \begin{equation}
    \label{aptau}
    \int_{\Rz^d \times (\tau,T-2\tau)} \left( |\partial_t  u_\tau|^2 + |\nabla  u_\tau|^2 + | u_\tau|^p\right) dx\, dt \leq C
  \end{equation}
where $C$ depends on $u^0$, $u^1$, and $T$ but not on $\epsi$ nor $\tau$.
\end{lemma}

\subsection{Passage to the limit as $\boldsymbol \tau \boldsymbol \to \boldsymbol 0$}\label{passage}

In order to conclude the proof of Lemma \ref{ee} we need to show that the time-discrete energy estimate \eqref{aptau} passes to the limit as $\tau \to 0$ (for fixed $\epsi>0$). To this aim, we check the $\Gamma$-convergence\cite{DeGiorgi-Franzoni75} $
  I_{\epsi\tau} \stackrel{\Gamma}{\to} I_\epsi$  with respect to the the weak topology of 
$$\VV\define H^1(0,T;L^2(\Rz^d))\cap L^2(0,T;H^1(\Rz^d))\cap L^p(\Rz^d\times (0,T))$$
where, clearly, the $(n+1)$-vector $(u_0,\dots,u_n)$ is intended to be identified with its piecewise affine interpolant $u_\tau$. More precisely, we prove the following.

\begin{lemma}[$\Gamma$-convergence]\label{g} Let 
  \begin{gather}
    \UU_{\text{\rm affine}} \define \{ u : [0,T]\to U \ : \ u \ \ \text{is piecewise affine on the time partition}\}\nonumber
  \end{gather}
and define the functionals $G_{\epsi},\, G_{\epsi\tau}: \UU \to [0,\infty]$ as
\begin{align}
&G_\epsi = I_\epsi \ \ \text{on} \ \ \KK(u^0,u^1) \ \ \text{and} \ \ G_{\epsi}\define \infty \ \ \text{elsewhere},\nonumber\\
&G_{\epsi\tau}(u_\tau)\define I_{\epsi\tau}(u(0),u(\tau), \dots,u(T)) \ \ \text{if} \ \ u_\tau \in \UU_{\text{\rm affine}}\cap \KK_\tau(u^0,u^1) \nonumber\\
&\ \ \text{and} \ \ G_{\epsi\tau}\define \infty \ \ \text{elsewhere}.\nonumber
\end{align}
Then,
\begin{gather}
  G_{\epsi\tau} \stackrel{\Gamma}{\to} G_\epsi \ \ \text{with respect to the the weak topology of} \ \VV.\nonumber
\end{gather}
\end{lemma}

The proof of this lemma is classically divided into proving the $\Gamma$-$\liminf$ inequality (Subsection \ref{gliminf}) and checking the existence of a recovery sequence (Subsection \ref{recov}). 

Before going on, let us remark that 
\begin{equation}
  \rho_\tau, \, \ove \rho_\tau, \, \ove \rho_\tau(\cdot -\tau), \, \ove \rho_\tau(\cdot -2\tau) \to \left( t \mapsto \e^{-t/\epsi}\right) \  \  \text{strongly in} \ \ L^\infty(0,T),\label{conv_rho}
\end{equation}
the convergence of $\rho_\tau$ being actually strong in $W^{1,\infty}(0,T)$.

\subsubsection{\bf $\boldsymbol \Gamma$-$\boldsymbol \liminf$ inequality}\label{gliminf}

Assume to be given $u_\tau \in \UU_{\text{\rm affine}}$ such that $u_\tau \to u$ with respect to the the weak topology of $ \VV$ and $\liminf_{\tau\to 0}G_{\epsi\tau}(u_\tau)<\infty$. 

 Let us denote by $\tilde u_\tau\in C^1([0,T];U)$ the piecewise-quadratic-in-time interpolant of $(u_0,\dots, u_n)\define (u_\tau(0), \dots, u_\tau(n\tau))$ defined by the relations
\begin{gather}\tilde u_\tau(t)\define
  u_\tau(t)  \ \  \text{for} \ \ t \in [0,\tau] \nonumber \\
\text{and} \ \ 
  \partial_{t} \tilde u_\tau(t) = \alpha_\tau(t)\partial_{t}  u_\tau(t) + (1 -\alpha_\tau(t))\partial_{t}  u_\tau(t-\tau) \ \ \text{for} \ \ t \in (\tau,T]
\end{gather}
where we have used the notation $ \alpha_\tau(t)\define  \alpha_i(t)$ for $t \in ((i-1)\tau,i\tau]$, $i=1,\dots,n$. We preliminarily observe that 
\begin{equation}\label{pre} \partial_t \tilde u_\tau(t) = \partial_t u_\tau(t - \tau) + \tau \alpha_\tau(t) \partial_{tt}\tilde u_\tau(t) \quad \forall t \in (\tau,T],
\end{equation}

As $\liminf_{\tau\to 0}G_{\epsi\tau}(u_\tau)<\infty$, by possibly extracting some not relabeled subsequence, we have that $u_\tau(0)=u^0$ and
$$\limsup_{\tau \to 0 }\left(\int_\tau^T \ove \rho_\tau \frac12 |\partial_{tt} \tilde u_\tau|^2 + \int_{\tau}^{T-2\tau} \frac{1}{\epsi^2}\ove \rho_\tau(\cdot +2\tau) \phi(\ove u_\tau) \right)< \infty.$$
Then, owing to convergence \eqref{conv_rho} and using
$$\delta u_k = u^1 + \sum_{i=2}^k \tau \delta^2 u_i \quad k=2,\dots, n,$$
we have that, for $\tau$ small,
$$\int_\tau^T  |\partial_{tt} \tilde u_\tau|^2 + \int_\tau^T  |\partial_{t} u_\tau|^2 + \int_{\tau}^{T-2\tau}  \phi(\ove u_\tau) \leq C $$
where $C$ depends on $u^1$, $T$, and $\epsi$. Hence, by possibly further extracting a not relabeled subsequence, we have that 
\begin{align}
\nabla \ove u_\tau \to \nabla u \quad &\text{weakly in} \ \ L^2(\Rz^d\times (0,T)),\label{conv_u1}\\
\ove u_\tau \to u \quad &\text{weakly in} \ \  L^p(\Rz^d\times (0,T)),\label{conv_u2}\\
 u_\tau \to u \quad &\text{weakly in} \ \ H^1(0,T;L^2(\Rz^d)),\label{conv_u3}\\
 \tilde u_\tau \to v\quad &\text{weakly in} \ \ H^2(0,T;L^2(\Rz^d)),\label{conv_u4}\\
\partial_t \tilde u_\tau(t) \to v_t(t) \quad &\text{weakly in} \ \ L^2(\Rz^d) \ \ \text{for all} \ \ t \in [0,T].\label{conv_u5}
\end{align}

We shall prove that, indeed, $v\equiv u$. To this aim, fix $w\in L^2(\Rz^d\times (0,T)) $ and compute that 
\begin{align}
& \int_0^T (\partial_t \tilde u_\tau - u_t,w)\stackrel{\eqref{pre}}{=}  \int_\tau^T (\partial u_\tau(\cdot - \tau) + \tau \alpha_\tau \partial_{tt}\tilde u_\tau - u_t,w) + \int_0^\tau (u^1 - u_t,w)\nonumber\\
&= \int_0^T (\partial_t u_\tau - u_t,w) - \int_\tau^T(\partial_t u_\tau -\partial_t u_\tau(\cdot - \tau),w ) + \tau \int_\tau^T\alpha_\tau (\partial_{tt}\tilde u_\tau,w) \nonumber\\
&=  \int_0^T (\partial_t u_\tau - u_t,w) - \tau \int_\tau^T(1-\alpha_\tau) (\partial_{tt}\tilde u_\tau,w) \to 0\nonumber  
\end{align}
where this convergence follows from \eqref{conv_u3}, $|\alpha_\tau|\leq 1$, and the boundedness of $\partial_{tt}\tilde u$. Namely, $\partial_t \tilde u_\tau \to u_t$ weakly in $L^2(\Rz^d\times (0,T))$ and  $v=u$.
In particular, owing to convergence \eqref{conv_u5} we have proved that $u^1= \partial_t \tilde u_\tau(0)=u_t(0) $ and $u\in \KK(u^0,u^1)$.

Eventually, we exploit Fatou's lemma and the convergences \eqref{conv_rho} and \eqref{conv_u4} in order to get that
\begin{align}
\int_0^T \e^{-t/\epsi}\frac12 | u_{tt}|^2 &\leq \liminf_{\tau \to 0} \int_\tau^T \ove \rho_\tau\frac12 |\partial_{tt} \tilde u_\tau|^2 = \liminf_{\tau \to 0} \sum_{i=2}^n \tau \rho_i \frac12 |\delta^2 u_i|^2,\nonumber\\
  \int_0^T \frac{1}{\epsi^2} \phi(u) &\leq \liminf_{\tau \to 0}\int_{2\tau}^{T-2\tau} \frac{1}{\epsi^2} \ove \rho_\tau(\cdot + 2\tau ) \phi(\ove u_\tau)= \liminf_{\tau \to 0} \sum_{i=2}^{n-2} \frac{\tau}{\epsi^2}\rho_{i+2} \phi(u_i).\nonumber
\end{align}
In particular, these last two inequalities ensure that 
\begin{align}
 G_\epsi(u) &\leq \liminf_{\tau \to 0}\left( \sum_{i=2}^n \tau \rho_i \frac12 |\delta^2 u_i|^2 + \sum_{i=2}^{n-2} \frac{\tau}{\epsi^2}\rho_{i+2} \phi(u_i)\right)\nonumber\\
&=  \liminf_{\tau \to 0} I_{\epsi\tau}(u_0, \dots,u_n) = \liminf_{\tau \to 0} G_{\epsi\tau}(u_\tau).\nonumber 
\end{align}

\subsubsection{\bf Existence of a recovery sequence}\label{recov}
Let us define the {\it backward mean operator}  
$$ M_\tau(u)(t)\define \media_{t - \tau}^{t} u(s)\, ds \quad \text{for} \ \ u \in \UU, \ t  >\tau .$$
Then, let $u\in \KK(u^0,u^1)$ be fixed and define $u_\tau$ by 
$$ u_0 \define u^0, \quad u_1 \define u^1 + \tau u^0, \quad u_i \define M_\tau(u)(i\tau) \quad \text{for} \ \ i=2, \dots,n.$$
We clearly have that $u_\tau \in \KK_\tau(u^0,u^1)$ and that both  $u_\tau$ and $ \ove u_\tau$ converge to $u$ strongly in $L^2(0,T;H^1(\Rz^d))\cap L^p(\Rz^d \times (0,T))$. Moreover, one can check that
\begin{align}
&\int_0^T |\partial_t u_\tau - u_t|^2 \nonumber\\
&= \int_0^\tau |u^1 - u_t|^2 +
 \sum_{i=2}^n \int_{(i-1)\tau}^{i\tau} \left| \frac{1}{\tau} M_\tau (u - u(\cdot - \tau)) - u_t\right|^2 \nonumber\\
&=\int_0^\tau |u^1 - u_t|^2 + \sum_{i=2}^n \int_{(i-1)\tau}^{i\tau} \left|\media_{(i-1)\tau}^{i\tau}M_\tau(u_t)(s)ds - u_t(t)\right|^2dt.\label{glimsup0}
\end{align}
Hence, as one has that $M_\tau(u_t)\to u_t$  strongly in $ L^2(\Rz^d \times (0,T))$,
we conclude that $u_\tau \to u $ strongly in $H^1(0,T;L^2(\Rz^d))$. In particular, we have checked that $u_\tau\to u $ weakly in $\VV$.

Next, we  compute that 
\begin{align}
&\sum_{i=2}^{n-2} \frac{\tau}{\epsi^2} \rho_{i+2} \phi(u_i) \nonumber\\
&= \sum_{i=2}^{n-2}\frac{1}{\epsi^2}{ \rho_{i+2} }
  \int_{(i-1)\tau}^{i\tau} \Big( \phi(\ove u_\tau(t)) - \phi(u(t))\Big) dt +  \sum_{i=2}^{n-2}\frac{1}{\epsi^2}{ \rho_{i+2} }
  \int_{(i-1)\tau}^{i\tau} \phi(u(t))dt\nonumber\\
&\leq  \sum_{i=2}^{n-2} \int_{(i-1)\tau}^{i\tau}\frac{ 1}{\epsi^2}\ove \rho_{\tau}(\cdot + 2\tau)  (A \ove u_\tau, \ove u_\tau - u)  + \sum_{i=2}^{n-2} \int_{(i-1)\tau}^{i\tau}\frac{ 1 }{\epsi^2}\ove \rho_{\tau}(\cdot + 2\tau) \phi(u)
\end{align}
where we have exploited the convexity of $\phi$.
In particular, by taking the $\limsup$ as $\tau \to 0$ and recalling that $\ove u_\tau \to u$ strongly in $L^2(0,T;H^1(\Rz^d))\cap L^p(\Rz^d \times (0,T))$ and \eqref{conv_rho}, we have that
\begin{equation}
  \label{glimsup}
 \limsup_{\tau \to 0} \left( \sum_{i=2}^{n-2} \frac{\tau}{\epsi^2} \rho_{i+2} \phi(u_i) \right) \leq \int_0^T \frac{1}{\epsi^2}\e^{-t/\epsi} \phi(u).
\end{equation}

We shall now specialize the argument of estimate \eqref{glimsup0} to the control of the difference of second-order derivatives in time by computing
\begin{align}
 & \sum_{i=2}^n \int_{(i-1)\tau}^{i\tau} \left| \delta^2 u_i - u_{tt}(t)\right|^2dt = \sum_{i=2}^n \int_{(i-1)\tau}^{i\tau} \left| \frac{u_i - 2u_{i-1} + u_{i-2}}{\tau^2} - u_{tt}(t)\right|^2dt  \nonumber\\
&= \sum_{i=2}^n \int_{(i-1)\tau}^{i\tau} \left| \frac{1}{\tau^3} \int_{(i-1)\tau}^{i\tau}(u - u(\cdot -\tau)) - \frac{1}{\tau^3}\int_{(i-2)\tau}^{(i-1)\tau}(u - u(\cdot - \tau)) -u_{tt}(t)\right|^2dt \nonumber\\
&= \sum_{i=2}^n \int_{(i-1)\tau}^{i\tau} \left| \frac{1}{\tau^3} \int_{(i-1)\tau}^{i\tau}\left(\int_{s-\tau}^s \Big( u_t - u_t(\cdot-\tau) \Big) \right) ds -u_{tt}(t)\right|^2dt \nonumber\\
&= \sum_{i=2}^n \int_{(i-1)\tau}^{i\tau} \left| \media_{(i-1)\tau}^{i\tau} \left( \media_{s-\tau}^sM_\tau(u_{tt})(r)dr\right)ds- u_{tt}(t)\right|^2dt \to 0\label{glimsup2}
\end{align}
where the convergence to $0$ is ensured by the fact that $M_\tau(u_{tt}) \to u_{tt} $ strongly in $ L^2(\Rz^d \times (0,T))$. In particular, the convergence in \eqref{glimsup2} amounts to say that (see Subsection \ref{gliminf} for the definition of $\tilde u_\tau$)
$$\partial_{tt} \tilde u_\tau \to  u_{tt} \quad \text{strongly in} \ \ L^2(\Rz^d\times (0,T)).$$

Finally, by collecting \eqref{glimsup}-\eqref{glimsup2} we have proved that
\begin{align}
&\limsup_{\tau \to 0}G_{\epsi\tau}(u_\tau) = \limsup_{\tau \to 0}I_{\epsi\tau}(u_0, \dots, u_n) \nonumber\\
&= \limsup_{\tau \to 0} \left( \sum_{i=2}^n \tau \rho_{i} \frac12 |\delta^2u_i|^2 +  \sum_{i=2}^{n-2} \frac{\tau}{\epsi^2} \rho_{i+2} \phi(u^i)\right)\nonumber\\
&= \limsup_{\tau \to 0} \left( \int_{2\tau}^T\ove \rho_\tau \frac12 |\partial_{tt}\tilde u_\tau|^2 + \int_{2\tau}^{T-2\tau}\frac{1}{\epsi^2}{\ove \rho_\tau(\cdot +2\tau)} \phi(\ove u_\tau)\right)\nonumber\\
& \leq \int_0^T \e^{-t/\epsi}\left( \frac12 |u_{tt}|^2 + \frac{1}{\epsi^2}\phi(u) \right) = G_{\epsi}(u).\nonumber
\end{align}
Namely, $u_\tau$ is a recovery sequence for $u$.


\subsubsection{\bf Proof of Lemma \ref{ee}}\label{lemmaee}
This is just an easy application of the Fundamental Theorem of $\Gamma$-convergence.\cite{DalMaso93} The minimizers $u_\tau$ of $G_{\epsi\tau}$ fulfill estimate \eqref{aptau} and are hence weakly precompact in $\VV$. As $G_{\epsi\tau}\stackrel{\Gamma}{\to} G_\epsi$ with respect to the same topology by Lemma \ref{g}, we have that $u_\tau \to u$ weakly in  $\VV$ where $u$ is the unique minimizer of $G_\epsi$. Finally, the bound \eqref{aptau} passes to the limit and the energy estimate \eqref{ap} holds.


\section{Proof of the Conjecture}\label{limit}
\setcounter{equation}{0}

Let us now come to the proof of Theorem \ref{main}. Owing to the energy estimate \eqref{ap} and by possibly implementing a diagonal extraction argument in combination with the {\sc Aubin-Lions} Lemma we obtain, for some not relabeled subsequence,
\begin{align}
  u^\epsi \to u \quad &\text{weakly in}  \ \ \VV \ \ \text{and strongly in} \ \ C([0,T];L^2(B_k)) \ \ \text{for all} \ \ k >0\label{conv1}
\end{align}
where $B_k\define \{ x \in \Rz^d \ : \ |x|\leq k\}$. Hence, $u^\epsi \to u $  and $|u^\epsi|^{p-2}u^\epsi \to |u|^{p-2}u $ pointwise almost everywhere. In particular, again from the energy estimate \eqref{ap} we have that \cite[Prop. 3.10, p. 272]{Visintin96}
\begin{gather}
|u^\epsi|^{p-2}u^\epsi \to |u|^{p-2}u \quad \text{strongly in} \ \ L^{q}(B_k \times (0,T))\nonumber\\
\text{for all} \ \ k >0\ \ \text{and} \ \   q \in [1,p').\label{conv3}
\end{gather}

We shall now check that $u$ solves \eqref{PDE} in the distributional sense. To this aim, let $w \in C^\infty_{0}(\Rz^d \times (0,T))$ and define $v\define \e^{t/\epsi}w$. We have that
$$ v_{tt} = \e^{t/\epsi} w_{tt}+ \frac{2}{\epsi}\e^{t/\epsi}w_t + \frac{1}{\epsi^2} \e^{t/\epsi} w.$$
As $v \in  C^\infty_{0}(\Rz^d \times (0,T)) \subset \KK(0,0)$, from the variational equality \eqref{var} one obtains
\begin{align}
 &0= \int_0^T \e^{-t/\epsi}\left( (u^\epsi_{tt} ,v_{tt} )  +\frac{1}{\epsi^2} (Au^\epsi ,v )\right) \nonumber\\
&= \int_0^T \left( (u^\epsi_{tt},w_{tt})  + \frac{2}{\epsi}(u^\epsi_{tt},w_t) + \frac{1}{\epsi^2} (u^\epsi_{tt},w) +\frac{1}{\epsi^2} (A u^\epsi,w)\right).\nonumber
\end{align}
Hence, we compute that 
\begin{align}
 &\int_0^T \Big(  (u^\epsi,w_{tt}) + (A u^\epsi,w)\Big) = \epsi^2 \int_0^T (u^\epsi_t,w_{ttt})  + 2\epsi \int_0^T (u^\epsi_t,w_{tt})\nonumber\\
&=\int_0^T(u^\epsi_t,\epsi^2 w_{ttt} + 2\epsi w_{tt}) .\nonumber
\end{align}
By passing to the limit as $\epsi \to 0$ and using the convergences \eqref{conv1}-\eqref{conv3} and the compactness of the support of $w$ we have that
$$\int_0^T \big(  (u,w_{tt}) + (A u,w)\big)=0.$$
Namely, $u$ solves \eqref{PDE} in the distributional sense.

\subsection{Comments and extensions}\label{comments}

As the limit function $u$ belongs to the space $\VV$, we have that  equation \eqref{PDE} is solved in the following stronger sense
  \begin{align}
    \int_{\Rz^d \times (0,T)}\left(-  u^\epsi_{t} w_{t} + \nabla u^\epsi\cdot\nabla w + \frac{p}{2}|u|^{p-2}u\, w\right) =0 \quad
 \text{for all} \ \ w \in \VV.
  \end{align}

Moreover, for $p<2^*$ Problem \eqref{PDE} admits a unique strong solution. In this case, the convergence $u^\epsi \to u$ clearly holds for the whole sequence.

By inspecting the proof of Theorem \ref{main} one realizes that the argument follows by assuming the weaker regularity on initial data $$
u^0, \, u^1 \in U.$$ 
Moreover, the choice of $\Rz^d$ in the definition of $I_\epsi$ can be replaced by some $\Omega \subset \Rz^d$ (possibly along with Dirichlet or even mixed boundary conditions) with no particular intricacy.

The power nonlinearity in $I_\epsi$ can be replaced by a more general convex differentiable function $F$ of polynomial growth. In particular, we can ask for
$$  \frac{1}{C}|r|^p \leq F(r) + C, \quad  |F'(r)|^{p'}\leq C (1 +  |r|^{p}) \quad \forall r\in \Rz.$$
The coercivity requirement can be omitted if $F$ is of subcritical growth, namely $p< 2^*$. An additional quadratic term in $F$ can also be considered. In particular, the Conjecture  \ref{1.1} holds true in the case of semilinear Klein-Gordon equations as well.

\section*{Acknowledgment}
Comments by Matthias Liero, Enrico Serra, and Paolo Tilli on an earlier version of the paper are gratefully acknowledged. Partial support from FP7-IDEAS-ERC-StG Grant \#200497 (BioSMA) is also acknowledged.


\def\cprime{$'$} \def\cprime{$'$}

\end{document}